\documentclass[journal]{IEEEtran}
\usepackage{graphicx}
\ifCLASSINFOpdf
\else
\fi
\usepackage{color}

\begin{document}
\title{Graph Partitioning with Fujitsu Digital Annealer}

\author{
	\IEEEauthorblockN{Yu-Ting~Kao\IEEEauthorrefmark{1}, Hsiu–Chuan~Hsu\IEEEauthorrefmark{1}\IEEEauthorrefmark{2}\thanks{email: hcjhsu@nccu.edu.tw}
		\\\IEEEauthorrefmark{1} Graduate Institute of Applied Physics, National Chengchi University, Taipei, Taiwan
	\\\IEEEauthorrefmark{2}Department of Computer Science, National Chengchi University, Taipei, Taiwan 
}

\thanks{}% <-this % stops a space
\thanks{}% <-this % stops a space
\thanks{}}


\markboth{}%
{}



\maketitle

\begin{abstract}
Graph partitioning, or community detection, is the cornerstone of many fields, such as logistics, transportation and smart power grids. Efficient computation and efficacious evaluation of communities are both essential, especially in commercial and industrial settings. However, the solution space of graph partitioning increases drastically with the number of vertices and subgroups. With an eye to solving large scale graph partitioning and other optimization problems within a short period of time, the Digital Annealer (DA), a specialized CMOS hardware also featuring improved algorithms, has been devised by Fujitsu Ltd. This study gauges Fujitsu DA's performance and running times. The modularity was implemented as both the objective function and metric for the solutions. The graph partitioning problems were formatted into Quadratic Unconstrained Binary Optimization
(QUBO) structures so that they could be adequately imported into the DA. The DA yielded the highest modularity among other studies when partitioning Karate Club, Les Miserables, American Football, and Dolphin. Moreover, the DA was able to partition the Case 1354pegase power grid network into 45 subgroups, calling for 60,930 binary variables, whilst delivering optimal modularity results within a solving time of roughly 80 seconds. Our results suggest that the Fujitsu DA can be applied for rapid and efficient optimization for graph partitioning.
\end{abstract}


% \fi
\IEEEpeerreviewmaketitle

\section{Introduction} \label{intro}
\IEEEPARstart{T}{he} industry and the academic sector are permeated with optimization problems \cite{Liao2023}. For business decisions, it is essential to minimize cost while maximizing profit, and for physics, it is often required to find the solution that minimizes the energy given a Hamiltonian that describes the system. Furthermore, optimization problems are evidenced in everyday-life scenarios, even to the extent of minute or mundane tasks such as putting together work schedules for employees. Although certain problems can be easily solved with analytical methods or classical computing, many are NP-hard, which means solving such problems within polynomial time are seldom encountered for classical computers and that the time complexity of such problems scale exponentially. 

Simulated Annealing (SA) \cite{SA} is a popular method among a myriad of other approaches attempting to tackle this issue. In quantum physics, SA has been applied to solve for the energy states for quantum Ising models, known to be an NP problem for classical computers \cite{Fu1986}. Outside the realm of quantum physics, it has been shown that many combinatorial optimization problems can be mapped to Ising models \cite{Lucas2014}\cite{Glover2019} and solved by SA. Fujitsu's Digital Annealer (DA) \cite{DA1}\cite{DA2} is an application-specific CMOS hardware that performs an improved version of SA to solve for Ising-type problems. The specific hardware design enables the speed-up of the processing times. The DA possesses three major differences that discern itself from traditional SA. The DA initializes all runs from the same arbitrary state, eliminating the need to repeatedly calculate initial energies. It also features a parallel-trial scheme which executes all flips individually in parallel at every Monte-Carlo step, increasing the acceptance probability. Moreover, the DA employs dynamic offset, which is an escape mechanism to shorten time spent in local minima. Consequently, the DA can be regarded as an adequate apparatus when we need to solve exceptionally large, complex optimization problems in a limited time frame.

Our research is focused on using the Fujitsu DA to solve graph partitioning problems, also known as community detection. Community detection is pivotal when it concerns virtual microgrid detection in power distribution networks. In order to integrate renewable energy sources to the power grid, the centralized power distribution network currently implemented must be revamped to reduce power loss during transmission and attain higher efficiency for the accommodation of green energy \cite{GRID1}\cite{GRID2}\cite{Morstyn2023}. A number of classical, non-quantum-inspired algorithms have previously been proposed for community detection. They can be categorized as agglomerative or divisive \cite{Scott2014}\cite{newman_finding_commun...}. Agglomerative methods search for vertex pairs that have similar characteristics and append the corresponding edges to an empty data structure, commencing from the most similar pair of vertices. Divisive methods, on the other hand, find dissimilar vertex pairs and remove their edges, starting from the least similar pair \cite{newman_finding_commun...}. One example of an agglomerative approach is Self-Avoiding Random Walk \cite{DEGUZZIBAGNATO20181046}. In a Self-Avoiding Random Walk, the walker is prohibited from revisiting a vertex in the same walk. The walker would constantly look for a viable path through vertices not yet traversed, and this single walk terminates when outstanding viable paths cease to exist. The number of steps in a path would be the total number of vertices and multiple walks are conducted. Edges connecting dissimilar vertex pairs are thought to be "bridges" or "highways" that channel two distinct communities, and are traversed more times than other edges within a certain community. An instance of a divisive approach is to remove the edges between the most dissimilar vertices according to the "betweenness " of a pair of vertices, introduced by Newman and Girvan \cite{newman_finding_commun...}. This algorithm looks for edges having the most "betweenness", an indicator that measures the contribution of an edge to the connection of communities. The betweenness can be the shortest-paths, random-walk betweenness or the current-flow betweenness, depending on the physical characteristics of the graph. %Breadth-first Search is utilized to find the shortest-path from a vertex to all the other vertices and betweenness is counted in the process.

Although these classical algorithms are effective, their time complexities are rather significant. In the aforementioned examples, the worst-case time complexities of $O(n^3log(n))$ and $O(n^3)$ are estimated for Self-Avoiding Random Walk \cite{DEGUZZIBAGNATO20181046} and the divisive method \cite{newman_finding_commun...} respectively, where $n$ is the total number of nodes or vertices in a graph. The DA aims to condense the required time to solve a community detection problem and attempts to secure maximum optimization. In this paper, we expand on our previous efforts\cite{kao2023solving} and place our emphasis on modularity optimization for graph partitioning and the influences of running time limits. In Section \ref{math}, we recapitulate and elaborate more in detail the mathematical formulations. In Section \ref{result}, we analyze in depth the results of classic graphs and examples of power distribution systems. Finally, the conclusion of this study is presented in Section \ref{concl}.

\section{Formulation} \label{math}
Several combinatorial optimization problems can be mapped to the Ising Model and the optimal configuration may be obtained by minimizing the model's energy, as in (\ref{eq:Ising}):
\begin{equation} \label{eq:Ising}
E(S) = -\sum_{i,j}w_{i,j}s_is_j-\sum_ih_is_i
\end{equation}
where the binary spin variables $s_i$ and $s_j$ are either $+1$ (up) or $-1$ (down), $w_{i,j}$ is the exchange energy or interaction between spins at node $i$ and $j$, $h_i$ is the external magnetic field, and $S$ is the set of all the binary variables. Should we require the utilization of binary variables $b_i\in{0,1}$, we can transform the Ising representation into a Quadratic Unconstrained Binary Optimization (QUBO) format by the conversion $b_i = \frac{s_i + 1}{2}$. The mathematical formulations in our work are expressed in QUBO form as required by the Fujitsu DA input structure.

%The graph partitioning problem can be mapped to the Ising model since we are solving for the best configuration of each node in a graph. 
Graph partitioning can be achieved by optimizing modularity, which can be mapped to QUBO form. Below, we explain the formulation in detail. 
The modularity $Q$ as the metric for quantifying the quality of communities is \cite{newman_finding_commun...}\cite{Blondel2008}\cite{Arenas_2008}:
\begin{equation} \label{eq:mod}
Q = \frac{1}{2m}\sum_{ij}\left(A_{ij} - \gamma\frac{k_ik_j}{2m} \right)\delta(c_{i},c_{j})\
\end{equation}
where $m = \frac{1}{2}\sum_{ij}^{}A_{ij}$ is the total number
of edges accounting for edge weights, $A_{ij}$is the coefficient in row $i$ and column $j$ of the graph's adjacency matrix, $\gamma$ is the resolution parameter which is set to 1 in our work, $k_{i}$ is the degree of node $i$ accounting for edge weights, and $\delta(c_{i},c_{j})=1$ if node $i$ and node $j$ are allocated to the same community otherwise 0.

The concept of modularity was coined by Newman and Girvan \cite{newman_finding_commun...}. It's value ranges from 0 to 1 and the closer the modularity is to 1 the better the partition is for the graph. Modularity essentially compares the proportion of edges in the same community to the probability that the edges in the graph are connected based on randomness.

The modularity can be utilized as a metric and be calculated after
partitioning a graph by any method for the evaluation of the resultant community structure. In this study, we use
modularity as an objective function in the QUBO form for
Fujitsu DA to conduct modularity optimization. To be more specific, the optimization tasks were performed on the third generation of the Fujitsu DA.

Given a graph $G$ with $n$ nodes and partitioned into $K$ groups,
we let
$\hat{x_{i}} = (x_{i0},\ x_{i1},x_{i2},\ldots,x_{ik},\ldots x_{i(K - 1)})$
be the vector associated with node $i$ and let $x_{ik}$ be a binary
variable, where $K$ is the total number of communities partitioned. If node $i$ is in group $k$, $x_{ik} = 1$ and all the
other binary variables in vector $\hat{x_{i}}$ would be zero. Since Fujitsu DA minimizes the energy, our objective function would be:
\begin{equation} \label{eq:Matrix}
M = -\textbf{ x}^{T}Q\textbf{x}
\end{equation}
where the vector
$\textbf{x} = (x_{00},x_{01},x_{02},\ldots,x_{ik},\ldots x_{n - 1,K - 1})^T$, to
maximize modularity. The index $i$ denotes the vertex and $k$ denotes the group. In this formalism, the number of binary variables required is $nK$. The size of all possible configurations, $2^{nK}$, scales exponentially with the number of binary variables.

In addition to the objective function, two constraints are needed for reasonable partitions. We require each node to be allocated to only one group and each of our groups to have at least one node to avoid empty groups. The two constraints are respectively:
\begin{equation} \label{eq:c1}
\sum_{k = 0}^{K - 1}x_{ik} = 1
\end{equation}
and
\begin{equation} \label{eq:c2}
\sum_{i = 0}^{n - 1}x_{ik} \geq 1.
\end{equation}
The first constraint, written in QUBO form becomes:
\begin{equation} \label{eq:c1-2}
C_1 = \sum_{i = 0}^{n}\left( \left( \sum_{k = 0}^{K - 1}x_{ik} \right) - 1 \right)^{2}
\end{equation}
, while the second constraint, by the slack variable method \cite{osti_2004640}, is converted into:
\begin{equation} \label{eq:c2-2}
\sum_{i = 0}^{n - 1}x_{ik} = d + 1
\end{equation}
, where $d$ is an integer variable and $0 \leq d \leq n - 1$. To represent $d$ in numerical form, we use one-hot encoding and
express $d$ using $n$ binary variables. Hence, in QUBO form, the
constraint is expressed as:
\begin{equation} \label{eq:c2-3}
C_2 = \sum_{k = 0}^{K - 1}\left( \left( \sum_{i = 0}^{n - 1}x_{ik} \right) - \left( \sum_{i = 0}^{n - 1}ix_{ik} \right) - 1 \right)^{2}.
\end{equation}
The full Hamiltonian that is imported into Fujitsu DA is:
\begin{equation} \label{eq:H}
H = M + \lambda_{c_1}C_1 + \lambda_{c_2}C_2
\end{equation}
, where $\lambda_{c_1}$ and $\lambda_{c_2}$ are the penalty multipliers
of constraints $C_1$ and $C_2$ respectively and are set empirically
depending on the problem \cite{GRID2}. Fujitsu DA provides a one-way one-hot constraint interface as well as an inequality constraint separation feature that assist users in tweaking these multipliers, which could reduce time spent on experimenting with various parameter values. One of the important parameters that control the solution quality is the maximum running time for the DA, controlled by the parameter $time\_limit\_sec$ in the Application Programming Interface (API). A longer running time permits the solver to visit more configurations and it would be more likely for it to locate the global minima. The actual solver running time is returned under the parameter $solve\_time$. It was found that there was a difference of a few seconds in those two parameters; thus, both values are reported in the Section \ref{result}.

\section{Results} \label{result}
\subsection{Test Graphs}
We compared the highest modularity obtained from Fujitsu DA with those reported in two other publications\cite{newman_finding_commun...}\cite{DEGUZZIBAGNATO20181046}. The results from these two publications were unweighted graphs.
Experiments were conducted with four graphs: Karate Club (KC), Les Miserables (LM), American Football (AF), and Dolphin (D). Karate Club and Les Miserables can either be weighted graphs or unweighted graphs, whilst others are unweighted graphs. Karate Club depicts a social interaction network in a karate club at an American university \cite{zachary_karate}, Les Miserables is a network of character co-appearances in the novel Les Miserables \cite{knuth1993stanford}, American Football describes a network of football games among division IA universities in the US \cite{football_newman}, and Dolphin exhibits the social network of bottlenose dolphins \cite{lusseau2003bottlenose}. Karate Club, Les Miserables, and American Football graphs were all imported from the Python library NetworkX version 3.1 \cite{SciPyProceedings_11}, and the data for Dolphin was downloaded from the website Network Data Repository \cite{nr}.  We provide both weighted and unweighted graph modularity results where it is applicable. We set the parameter $"time\_limit\_sec"$ to 10 for each graph executed, which would limit the running time to approximately 10 seconds for each graph. We also partitioned each graph with the parameter $"time\_limit\_sec"$ set to 80, and verified that the results were identical to their 10-second counterparts.

The optimal modularity and the corresponding number of partitions are shown in Table \ref{tab:graph_res}. Fujitsu DA performed outstandingly, yielding the most optimal results for all four graphs and satisfying both constraints. KC results for a weighted graph were reported in our previous publication \cite{kao2023solving} and we performed additional experiments on KC, LM, AF, and D for this paper.

\begin{table}[htb]
    \centering
\caption{Network Test Results}
\label{tab:graph_res}
    \begin{tabular}{|c|c|c|c|c|} \hline 
         Q/K&  DA Weighted&DA Unweighted&  \cite{newman_finding_commun...}&  \cite{DEGUZZIBAGNATO20181046}\\ \hline 
         KC&  0.4449/4&\textbf{0.4198/4}&  0.4/2
&  0.4197/4\\ \hline 
         LM&  0.5667/6&\textbf{0.5600/6}&  0.54/11
&  0.5467/7
\\ \hline 
         AF&  N/A&\textbf{0.6046/10}&  N/A
&  0.6044/10
\\ \hline 
         D&  N/A& \textbf{0.5285/5}&  0.52/5
&  0.5277/5
\\ \hline
    \end{tabular}   
\end{table}

\subsection{Practical Examples}
As a practical application of our study, we conducted graph partitioning with Fujitsu DA on two power grid graphs in addition to the IEEE case studies provided in our previous work \cite{kao2023solving}. The edge weights were modelled as the inverse of the absolute value of impedance $\frac{1}{|r+jx|}$ \cite{inverse_of_impedance}, where $r$ was the resistance, $x$ was the reactance, and $j$ was $\sqrt{-1}$ here. The first power grid featured a distribution network of Taiwan Power Company (TPC) \cite{1208393}, with a total of 94 nodes and 96 edges. The second power grid was Case 1354pegase, which represented a portion of the European high voltage transmission network \cite{josz2016ac}\cite{Fliscounakis2013ContingencyRW}, with 1354 nodes and 1710 edges. We utilized the $nx.Graph()$ data structure of NetworkX version 3.1 \cite{SciPyProceedings_11} to construct a graph from the raw data \cite{1208393}, whilst we imported Case 1354pegase from the Python package PandaPower version 2.13.1 \cite{pandapower.2018}. We omitted parallel edges in our analyses for all of the graphs.

\subsubsection{Taiwan Power Company Distribution Network (TPC)}
\begin{figure}[tb]
    \centering
    \includegraphics[width=1.0\linewidth]{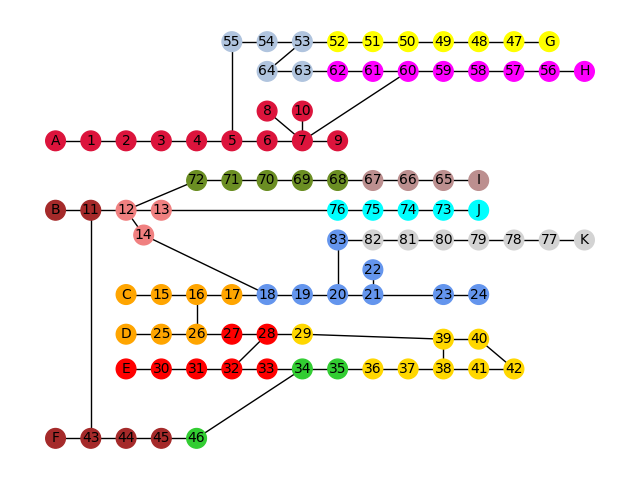}
    \caption{Partitioned TPC graph with the highest modularity 0.858432, with group number $K=15$, by Fujitsu Digital Annealer.}
    \label{fig:TAI_graph}
\end{figure}

\begin{figure}[htb]
    \centering
    \includegraphics[width=1.0\linewidth]{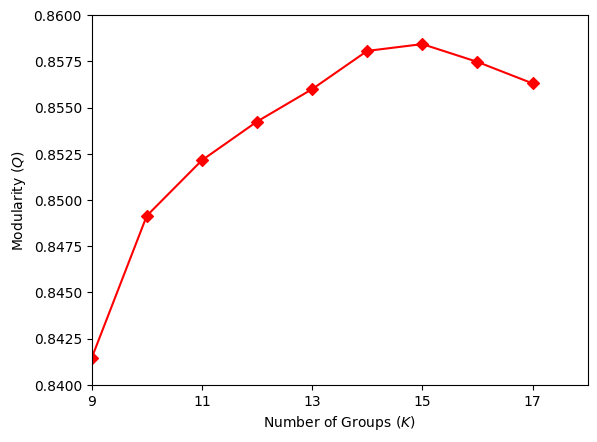}
    \caption{ Modularity versus number of groups for TPC graph solved with Fujitsu DA. The highest modularity is 0.858432 for 15 groups.}
    \label{fig:TAI_plot}
\end{figure}

The partitioned graph bearing the highest modularity for TPC is shown in Fig. \ref{fig:TAI_graph}, and the modularity against number of communities for TPC is plotted in Fig. \ref{fig:TAI_plot}. The highest modularity for TPC is 0.858432 for 15 groups. We set the parameter $"time\_limit\_sec"$ to 10 when partitioning TPC, which would limit the running time to approximately 10 seconds. We also partitioned TPC  with the parameter $"time\_limit\_sec"$ set to 80, and verified that the results were identical to its 10-second counterpart. 

\subsubsection{Case 1354pegase}
\begin{figure}[tb]
    \centering
    \includegraphics[width=1.0\linewidth]{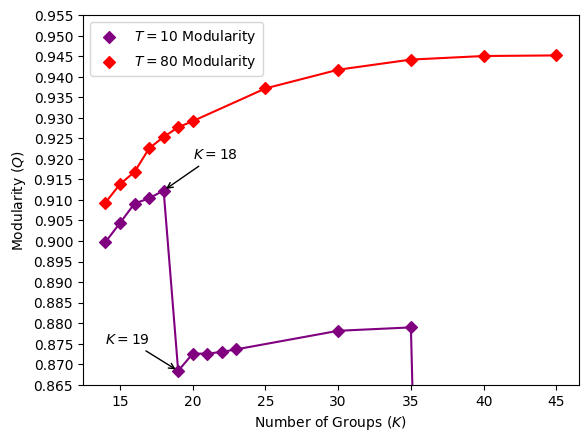}
    \caption{ Modularity versus number of groups for parameter $time\_limit\_sec\, (T)=10$ and $80$, solved with Fujitsu DA for Case 1354pegase. The modularity for 45 groups is 0.945187 at $time\_limit\_sec\, (T)=80$.}
    \label{fig:1354_plot}
\end{figure}

The modularity against number of communities for Case 1354pegase is shown in Fig. \ref{fig:1354_plot}. We were able to partition the graph up to 45 groups and the modularity value we obtained at 45 groups was 0.945187 when the parameter $time\_limit\_sec$ was set to 80. We initially ran Fujitsu DA with $time\_limit\_sec$ at 10, commensurate with TPC and previous test graphs, to study the impact of a larger graph on Fujitsu DA's performance. We subsequently observed a sharp drop in modularity performance at 19 and 36 groups, as in Fig. \ref{fig:1354_plot}, which indicated a requisite for longer execution times. For aesthetic purposes, we omitted results above 35 groups where $time\_limit\_sec$ was only 10.

\begin{figure}[tb]
    \centering
    \includegraphics[width=1\linewidth]{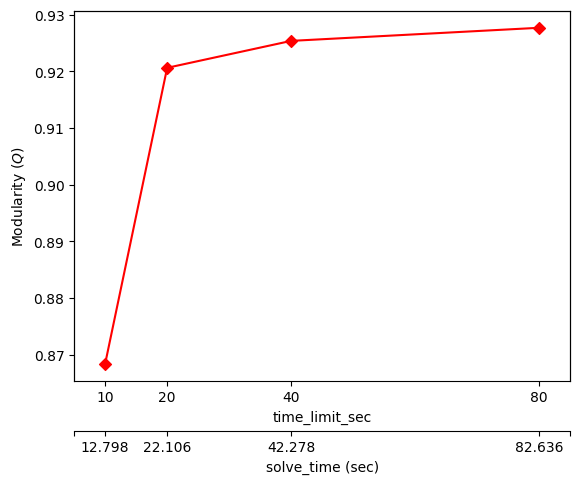}
    \caption{Modularity versus time limit and the solve time for 19 groups, solved with Fujitsu DA for Case 1354pegase. The modularity at $time\_limit\_sec=10$ is 0.868371, the modularity at $time\_limit\_sec=20$ is 0.92062, and the modularity at $time\_limit\_sec=80$ is 0.927689.  }
    \label{fig:1354 k19 time}
\end{figure}

\begin{figure}[tb]
    \centering
    \includegraphics[width=1\linewidth]{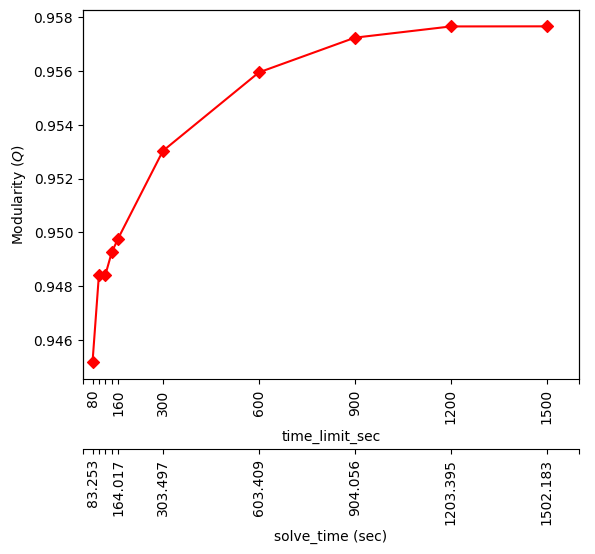}
    \caption{Modularity versus the time limit and the solve time for 45 groups, solved with Fujitsu DA for Case 1354pegase. The modularity at $time\_limit\_sec=80$ is 0.945187 and the modularity at $time\_limit\_sec=1500$ is 0.95764.}
    \label{fig:1354 k45 time}
\end{figure}

We further studied the relation between modularity and running time for 19 groups and 45 groups. 19-groups was where the drop occurred and deserved an examination, whilst 45-groups was thoroughly scrutinized in an attempt to find its maximum modularity. As aforementioned, the parameter $time\_limit\_sec$ limits the solving time; nonetheless, the actual solving time reported by Fujitsu DA in the entry called $solve\_time$ would exceed this limit by a few seconds. Hence, the actual solving times are also revealed in Fig. \ref{fig:1354 k19 time} and Fig. \ref{fig:1354 k45 time}. Despite this discrepancy, the actual solving times didn't deviate much from what was set for $time\_limit\_sec$, which is why we believe $time\_limit\_sec$ does in effect limit the running time and is generally indicative of the DA's actual execution time.


In 19-groups, the modularity underwent a stark escalation even though $time\_limit\_sec$ only increased by 10 seconds from 10 to 20, as shown in Fig. \ref{fig:1354 k19 time}. The modularity rose by a total of 6.8\% from 0.868371 at $time\_limit\_sec=10$ to 0.927689 at $time\_limit\_sec=80$. Once $time\_limit\_sec$ reached 20 or above, the increase in modularity became subtle, achieving only a 0.77\% growth in modularity from 0.92062 at $time\_limit\_sec=20$ to 0.927689 at $time\_limit\_sec=80$.

For 45-groups, we perceived that the modularity gradually converged to 0.95764, reaching 0.957637 after running for around 1200 seconds and attaining a maximum value of 0.95764 after approximately 1500 seconds of running time, as illustrated in Fig. \ref{fig:1354 k45 time}. In spite of the somewhat magnified depiction in Fig. \ref{fig:1354 k45 time}, the modularity merely increased 1.32\% from 0.945187 at $time\_limit\_sec=80$ to 0.95764 at $time\_limit\_sec=1500$. Therefore, we believe that the DA can still deliver desirable results in a limited time duration as short as approximately 80 seconds for a 45-group partition of the Case 1354pegase graph. A 45-group partition of the Case 1354pegase graph required $1354\times45=60,930$ binary variables in total, roughly half of the maximum 100,000 binary variables available on the third generation of the Fujitsu DA. Despite the large size, the DA was able to return optimal solutions within reasonable times.

\section{Conclusion} \label{concl}
We performed a thorough analysis regarding the ability of Fujitsu's Digital Annealer (DA) to explore the most optimal configurations in graph partitioning problems. The graph partitioning or community detection problem was transformed into the Quadratic Unconstrained Binary Optimization (QUBO) format and subsequently imported to Fujitsu DA. We evaluated the relation between modularity and the number of communities partitioned. We also shed light on the influence of running time on modularity optimization outcomes. Fujitsu DA performed exceptionally well in four test graphs. Fujitsu DA achieved a modularity of 0.4449 for Karate Club with weights, 0.5667 for Les Miserables with weights, 0.6046 for American Football, and 0.5285 for Dolphin. As practical power grid case studies, we evaluated a power distribution network from Taiwan Power Company (TPC) and also Case 1354pegase. The highest modularity obtained by Fujitsu DA was 0.858432 for 15 partitions of the TPC distribution network. For Case 1354pegase, a 45-group partition was accomplished with a modularity of 0.945187 within a running time of roughly 80 seconds. Further experiments revealed that the modularity attained at about 80 seconds differed little from 0.95764, the modularity obtained at about 1500 seconds of solving time, only a 1.32\% increase. Hence, we believe our findings in this paper can evince the potential of DA to rapidly solve community detection and power grid optimization problems among industries.

% or


\section*{Acknowledgments}
The authors would like to thank Prof. Yu–Cheng Lin and
Yu–Chen Shu for their insightful discussions. This work is
supported by National Science and Technology Council (NSTC) of Taiwan under the grant No. 112–2119–M–006–
004.

\ifCLASSOPTIONcaptionsoff
  \newpage
\fi



%\bibliographystyle{IEEEtran}
\bibliographystyle{ieeetr}
\bibliography{REVISEDmybib}

%\end{IEEEbiography}
\newpage
\begin{IEEEbiographynophoto}{Yu-Ting, Kao}
Yu-Ting, Kao achieved his B.S. in Electrophysics and B.S. in Money and Banking at National Chengchi University, Taipei City, Taiwan. He is currently a Research Assistant under Prof. Hsiu–Chuan, Hsu at the Graduate Institute of Applied Physics, National Chengchi University and is currently planning to pursue a Master's degree in science and technology.
\end{IEEEbiographynophoto}

%\newpage

\begin{IEEEbiographynophoto}{Hsiu–Chuan, Hsu}
Hsiu–Chuan, Hsu attained her Ph.D. in physics from Pennsylvania State University, USA. She is currently an Assistant Professor at the Graduate Institute of Applied Physics and Department of Computer Science at National Chengchi University, Taipei City, Taiwan.
\end{IEEEbiographynophoto}


%\vfill


\end{document}